\documentclass [11pt,oneside,a4paper,mathscr]{amsart}

\usepackage{amscd,amsmath,amssymb,euscript}
\usepackage[frame,cmtip,curve,arrow,matrix,line,graph]{xy}

\title{Hall algebras and curve-counting invariants}
\author{Tom Bridgeland}

\date{}

\newtheorem{thm}{Theorem}[section]

\newtheorem{cor}[thm]{Corollary}
\newtheorem{prop}[thm]{Proposition}
\newtheorem{lemma}[thm]{Lemma}
\newenvironment{pf}{\paragraph{Proof}}{\qed\par\medskip}
\theoremstyle{definition}
\newtheorem{defn}[thm]{Definition}

\newtheorem{thm*}[thm]{Theorem$^*$}

\newcommand{\comment}[1]{}

\renewcommand{\leq}{\leqslant}
\renewcommand{\geq}{\geqslant}
\newcommand{\isom}{\cong}
\newcommand{\tensor}{\otimes}
\newcommand{\onto}{\twoheadrightarrow}

\newcommand{\into}{\hookrightarrow}

\newcommand{\lRa}[1]{\xrightarrow{\ #1\ }}

\newcommand{\lra}{\longrightarrow}

\newcommand{\dual}{\vee}

\newcommand{\QQ}{\mathbb{Q}}

\newcommand{\D}{{D}}

\newcommand{\C}{\mathbb C}

\newcommand{\A}{\mathcal A}

\newcommand{\CC}{\mathcal C}
\newcommand{\OO}{\mathcal O}
\renewcommand{\P}{\mathcal{P}}

\newcommand{\LL}{\mathbb{L}}

\newcommand{\N}{\mathcal{N}}
\newcommand{\M}{\mathcal{M}}
\newcommand{\Q}{\mathcal{Q}}
\newcommand{\Z}{\mathbb{Z}}
\newcommand{\R}{\mathbf R}

\newcommand{\G}{\mathcal{G}}
\newcommand{\J}{\mathcal{J}}

\newcommand{\DD}{\mathbb{D}}

\newcommand{\lHom}{{\mathcal{H}om}}
\newcommand{\lExt}{{\mathcal{E}xt}}
\newcommand{\Coh}{\operatorname{Coh}}

\newcommand{\Hilb}{\operatorname{Hilb}}

\newcommand{\id}{\operatorname{id}}
\newcommand{\Ext}{\operatorname{Ext}}
\newcommand{\Hom}{\operatorname{Hom}}

\newcommand{\Spec}{\operatorname{Spec}}

\newcommand{\DT}{\operatorname{DT}}
\newcommand{\GL}{\operatorname{GL}}
\newcommand{\ch}{\operatorname{ch}}

\newcommand{\PT}{\operatorname{PT}}

\newcommand{\virt}{\operatorname{vir}}
\newcommand{\td}{\operatorname{td}}

\newcommand{\SF}{\operatorname{SF}}
\newcommand{\semi}{\operatorname{SS}}
\newcommand{\ad}{\operatorname{ad}}
\newcommand{\add}{\operatorname{Ad}}

\newcommand{\Hilbs}{\operatorname{Hilb}^\sharp}
\newcommand{\Hs}{\H^\sharp}
\newcommand{\sch}{\operatorname{Sch}}

\newcommand{\coker}{\operatorname{coker}}
\newcommand{\Knum}{N}

\newcommand{\Set}{\operatorname{Set}}
\newcommand{\op}{\operatorname{op}}
\renewcommand{\sc}{\operatorname{sc}}
\newcommand{\Qlast}{\Q\cap\CC}

\newcommand{\Fi}[1]{\operatorname{H}_\infty(#1)}

\newcommand{\kvar}[1]{K({\operatorname{Var/{#1}}})}

\newcommand{\kst}[1]{K({\operatorname{St/{#1}}})}
\newcommand{\RH}{\operatorname{H}}

\newcommand{\one}[1]{\operatorname{1}_{#1}}
\renewcommand{\O}[1]{\operatorname{1}^\OO_{#1}}
\newcommand{\RHreg}{\operatorname{H_{{reg}}}}
\newcommand{\RHsc}{\operatorname{H_{sc}}}

\begin{document}

\begin{abstract}
We use Joyce's theory of motivic Hall algebras to prove that
reduced Donaldson-Thomas  curve-counting invariants on Calabi-Yau
threefolds coincide
with stable pair  invariants, and that the  generating functions for these invariants
 are Laurent expansions of rational functions.
\end{abstract}
\maketitle

\section{Introduction}

In this paper we use Joyce's theory of motivic Hall algebras to
prove some basic properties of Donaldson-Thomas (DT)
curve-counting invariants on Calabi-Yau
threefolds. 
We prove that the reduced DT invariants coincide with the stable
pair invariants introduced by Pandharipande and Thomas, and that
the generating functions for these invariants are Laurent
expansions of rational functions. Similar results have been
obtained by Toda; we discuss the relationship with his work below.

The Hall algebra approach to  DT invariants relies  on a
fundamental result of Behrend. Recall \cite{Be} that Behrend
associates to any scheme $S$ of finite type over $\C$, a
 constructible function
\[\nu_S\colon S\to\Z,\]
with the property that if $S$ is a proper moduli scheme with a
symmetric obstruction theory, then the associated
virtual counting invariant
\[\#_{\virt}(S) := \int_{[S]^{\virt}} 1\]
coincides with the weighted Euler characteristic \[ \chi(S,\nu_S)
:= \sum_{n\in\Z} n\, \chi(\nu_S^{-1}(n)).\] This means that the
virtual count can be computed by motivic techniques involving
decomposing $S$ into a disjoint union of locally-closed
subschemes.

%

Suppose that $M$ is a smooth projective Calabi-Yau threefold over
$\C$. We include in this assumption the condition that
\[H^1(M,\OO_M)=0.\]
Let $N_1(M)$ denote the abelian group of curves in $M$ up to
numerical equivalence. It is a free abelian group of finite rank.
Set
\[N_{\leq 1}(M) = N_1(M) \oplus \Z.\] A coherent sheaf $E$ on $M$
supported in dimension $\leq 1$ has a Chern character
\[\ch(E)=(\ch_2(E),\ch_3(E))=(\beta,n)\in N_{\leq 1}(M).\]

Given a class $(\beta,n)\in N_{\leq 1}(M)$ there is an open and
closed subset of the Hilbert scheme
\[\Hilb_M(\beta,n)\subset
\Hilb_M\] parameterizing epimorphisms of coherent sheaves
\[f\colon\OO_M\onto E\] with $E$ of Chern character $(\beta,n)$. The
corresponding DT curve-counting invariant is
\[\DT(\beta,n)=\#_{\virt}\big(\Hilb_M(\beta,n)\big)\in \Z.\]
For a fixed curve class $\beta\in N_1(M)$ the scheme
$\Hilb_M(\beta,n)$ is known to be empty for $n\ll 0$, so the
generating function
\[\DT_\beta(q)=\sum_{n\in \Z} \DT(\beta,n) q^n\]
is a Laurent series.

In the case $\beta=0$ the scheme $\Hilb_M(\beta,n)$
parameterises zero-dimensional subschemes of $M$ of length $n$,
and the corresponding generating function $\DT_0(q)$ is known
\cite{BeF,LP,Li}. It is the Taylor expansion of the function
\[{M}(-q)^{\,\chi(M)},\] where $\chi(M)$ is the topological Euler
characteristic of $M$, and the MacMahon function
\[{M}(q)=\big.\prod_{k\geq 1}{ (1-q^k)^{-k}}\]  is the
generating function for three-dimensional partitions. Following
\cite{MNOP} one defines the reduced DT generating function to be
the quotient
\[\DT'_\beta(q)=\DT_\beta(q)/\DT_0(q).\]
It is a Laurent series  with integral coefficients.

 Pandharipande and Thomas  \cite{PT1} introduced new
curve-counting invariants by considering stable pairs. These are
maps of coherent sheaves \[f\colon \OO_M\to E\] such that $E$ is supported
in dimension one, and such that

\begin{itemize}
\item[(a)] $E$ has no zero-dimensional subsheaves,

\item[(b)] the cokernel of $f$ is zero-dimensional.
\end{itemize}
In fact, as we prove  in Lemma \ref{oneup} below, such maps can be
 though of as epimorphisms in the heart of a certain perverse
t-structure on the derived category  of $M$.

 It was proved in
\cite{PT1} that there is a fine moduli scheme parameterizing
stable pairs, which we will denote $\Hilb^\sharp_M$. For each
class $(\beta,n)\in N_{\leq 1}(M)$ there is an open and closed
subscheme
\[\Hilb^\sharp_M(\beta,n)\subset \Hilb^\sharp_M\]
parameterising stable pairs with $E$ of Chern character $(\beta,n)$. The
corresponding PT invariant is
\[\PT(\beta,n)=\#_{\virt}\big(\Hilb^\sharp_M(\beta,n)\big)\in \Z.\] Once
again, for fixed $\beta\in N_1(M)$  the generating function
\[\PT_\beta(q)=\sum_{n\in \Z} \PT(\beta,n) q^n\]
is a Laurent series. Note that in the case $\beta=0$ we have
$\PT_0\equiv 1$.

Our main result is as follows.

\begin{thm}
\label{conjdt} Let $M$ be a smooth projective Calabi-Yau threefold.
 Then for each  class $\beta\in
N_1(M)$
\begin{itemize}
\item[(a)] there is an equality of Laurent series
\[\DT'_\beta(q) = \PT_\beta(q).\]
\item[(b)] the series $\DT'_\beta(q)$ is the Laurent expansion of a
rational function in $q$, invariant under the transformation
$q\leftrightarrow q^{-1}$.
\end{itemize}
\end{thm}


%

 The first part of Theorem \ref{conjdt} (the DT/PT
correspondence) was conjectured by Pandharipande and Thomas
\cite[Conj. 3.3]{PT1}.
 The second part was conjectured earlier by Maulik, Nekrasov, Okounkov and Pandharipande \cite[Conj. 2]{MNOP}
  and is a crucial ingredient in the conjectural DT/GW correspondence
   \cite[Conj. 3]{MNOP}. 

There is another version of Theorem \ref{conjdt} involving
unweighted Euler characteristic invariants defined without a
Behrend function. This was first proved
 by Toda \cite{To1,To2}. Our argument
gives a second proof of this result, and 
Stoppa and Thomas \cite{ST} have given a third (and simpler) proof, which even extends to the
non-Calabi-Yau case.

 The approach to Theorem \ref{conjdt} we present here
owes a lot to Toda's work \cite{To1}, and we have used ideas from
his paper in several places. In fact Toda's argument could be
modified to give a proof of Theorem \ref{conjdt} if one could
establish a certain identity involving  Behrend weights for
objects of the derived category of $M$. The corresponding identity
for  coherent sheaves was established by Joyce and Song \cite{JS}
and plays a crucial role in this paper.

  The main difference between our approach and that of
\cite{To1} is that  Toda views
the scheme $\Hilb_M(\beta,n)$ as parameterizing stable rank 1 ideal sheaves
on $M$, whereas we consider it as parameterizing sheaves supported
in dimension $\leq 1$ together with the framing data of an
epimorphism from $\OO_M$. 

\smallskip
\noindent {\bf Acknowledgements} Thanks most of all to Dominic
Joyce who patiently explained many things about motivic Hall
algebras. In particular he pointed out how to apply his virtual
indecomposable technology  to replace the Kontsevich-Soibelman
integration map used in an earlier verion of this paper. I am also
very grateful to Yukinobu Toda who explained his approach to
Theorem \ref{conjdt} at the inaugural meeting of the HoRSe
seminar, and to Richard Thomas who provided a lot of encouragement
and technical know-how.  Finally, I'm  very happy to acknowledge
the entertaining and  useful conversations I've enjoyed with
Davesh Maulik, Rahul Pandharipande and Bal\'azs Szendr\H{o}i.

\smallskip
\noindent {\bf Notation} By a stack we mean an  algebraic (Artin)
stack. All schemes and stacks will be locally of finite type over
$\C$. Throughout $M$ will be a smooth complex projective threefold
with trivial canonical bundle and satisfying $H^1(M,\OO_M)=0$. We
write $\A=\Coh(M)$ for the abelian category of coherent sheaves on
$M$, and denote by $\M$ the  moduli stack of objects of $\A$. We
frequently use the same letter for an open substack of $\M$ and
the  corresponding full subcategory of $\A$. One other minor abuse
of notation: if $f\colon T\to S$ is a morphism of schemes, and $E$
is a sheaf on $S\times M$, we use the shorthand
 $f^*(E)$ for the pullback to $T\times M$, rather than the more correct $(f\times \id_M)^*(E)$.

\renewcommand{\H}{\mathcal{H}}


\section{Preliminaries}
We begin by assembling some basic facts that will be needed later.


\subsection{Curve classes}
We denote by $N_1(M)$ the abelian group of cycles of dimension 1
on $M$ modulo numerical equivalence. It is a free abelian group of
finite rank \cite[Ex. 19.1.4]{fulton}. An element $\beta\in
N_1(M)$ is effective, written $\beta\geq 0$, if it is zero or the
class of a curve.

\begin{lemma}
\label{eff} An element $\beta\in N_1(M)$ has only
finitely many decompositions of the form $\beta=\beta_1+\beta_2$ with  $\beta_i\geq 0$.
\end{lemma}

\begin{pf}
This follows immediately from \cite[Cor. 1.19]{KM}. \end{pf}

 We denote
by $\Knum(M)$ the numerical Grothendieck group of $M$. It is the
quotient of the Grothendieck group $K(M)$ by the kernel of the
Euler form. There is a subgroup
\[N_{\leq 1}(M)\subset N(M)\]
 generated by classes
of sheaves supported in dimension at most 1.

\begin{lemma}
\label{chuff}
The Chern character induces an isomorphism
\[\ch=(\ch_{2},\ch_3)\colon N_{\leq 1}(M) \lra N_1(M)\oplus \Z.\]
\end{lemma}

\begin{pf}
The Chern character defines a map to the Chow ring
\[\ch\colon K(M)\to A^*(M)\tensor\mathbb{Q}.\]
Let $K_{\leq 1}(M)\subset K(M)$ be the  subgroup spanned by
classes of sheaves  supported in dimension at most 1. The Chern
character of such a sheaf is equal to its Chern class for
dimension reasons, and in particular is integral. Thus for any
$\alpha\in K_{\leq 1}(M)$
\[\ch(\alpha)=(\ch_{2}(\alpha),\ch_{3}(\alpha))\in A_1(M)\oplus \Z.\]
By the Riemann-Roch theorem, the Euler form on $K(M)$ is given by
the formula
\[\chi(\alpha,\beta)=[\ch(\alpha^{\dual})\cdot\ch(\beta)\cdot\td(M)]_3.\]
It follows that an element $\alpha\in K_{\leq 1}(M)$ is  numerically trivial
precisely if $\ch_3(\alpha)=0$ and $\ch_{2}(\alpha)\cdot D=0$ for
any divisor $D$. Thus we get an injective map
\[\ch=(\ch_{2},\ch_3)\colon N_{\leq 1}(M)\to N_1(M)\oplus\Z.\]
Since $\OO_C(m)$
maps to $([C],\ch_3(\OO_C)+m)$, this map is also surjective.
\end{pf}

There is an effective cone
\[\Delta\subset N_{\leq 1}(M)\]
 consisting of classes of
sheaves. From now on we shall use the isomorphism  of Lemma
\ref{chuff} to identify $N_{\leq 1}(M)$ with the group
$N_1(M)\oplus \Z$. Under this identification the effective cone
becomes
\[\Delta= \big\{(\beta,n)\in N_1(M)\oplus \Z : \beta>0\text{ or
}\beta=0\text{ and } n\geq 0\big\}.\] Note also that the Euler
form $\chi(-,-)$ on $N(M)$ vanishes when restricted to $N_{\leq
1}(M)$ for dimension reasons. This trivial observation will be
very important in what follows.



\subsection{Stable pairs}
\label{pairs}

Define a full subcategory  \[\P=\Coh_0(M)\subset \A=\Coh(M)\]
consisting of sheaves supported in dimension 0. Let $\Q$ be the
full subcategory
\[\Q=\{ E\in \A : \Hom(P,E)=0 \text{ for all }P\in \P\}.\]
These subcategories  $(\P,\Q)$ form a torsion pair:
\begin{itemize}
\item[(a)] if $P\in\P$ and $Q\in\Q$ then $\Hom_\A(P,Q)=0$,
\item[(b)]every object $E\in\A$ fits into a  short exact
sequence
\[0\lra P\lra E\lra Q\lra 0\]
with $P\in \P$ and $Q\in\Q$.\end{itemize}
 Note that it follows from (a)
that the short exact sequence appearing in (b) is unique up to a
unique isomorphism.

We can define a t-structure on $\D(\A)$ by tilting the standard
t-structure with respect to this torsion pair \cite{HRS}. The resulting
t-structure has heart
\[\A^\sharp =\big\{E\in\D(\A): H^0(E)\in\Q,\  H^1(E)\in\P,\  H^i(E)=0\text{ for }i\notin\{0,1\}\big\}.\]
Note that $\Q=\A\cap\A^\sharp$ and, in particular, $\OO_M\in
\A^\sharp$.

\begin{lemma}
\label{oneup} \begin{itemize}
\item[(a)]If $f\colon \OO_M\to E$ is an
epimorphism in $\A^\sharp$ then   $E\in \Q$, and considering $f$ as a map in $\A$, one has $\coker(f)\in \P$.

\item[(b)] If $f\colon \OO_M\to E$
is a morphism in $\A$ with $E\in \Q$ and $\coker(f)\in\P$ then $f$ is an epimorphism in $\A^\sharp$.
\end{itemize}
\end{lemma}

\begin{pf}
If $f$ is an epimorphism in $\A^\sharp$ then there is a triangle
\begin{equation}
\label{tria}
K\lra \OO_M\lRa{f} E\lra K[1]\end{equation}
with $K\in\A^\sharp$.
Applying the long
exact sequence in cohomology gives \begin{equation}
\label{grimes}0\lra H^{0}(K) \lra \OO_M
\lRa{f} H^0(E)\lra H^1(K)\lra 0.\end{equation} It follows that $E=H^0(E)\in\Q$ and $\coker(f)=H^1(K)\in\P$.

   For (b) put the map $f$ in a triangle of the form \eqref{tria} and consider the associated
   long exact sequence
   \eqref{grimes}. Then $H^0(K)$ is a subsheaf of
$\OO_M$ and hence an object of $\Q$, and
\[H^1(K)=\coker(f)\in \P.\]
Therefore $K\in \A^\sharp$, the triangle is a short exact sequence in $\A^\sharp$ and hence $f$ is an epimorphism in $\A^\sharp$.
\end{pf}

Thus a stable pair in the sense of \cite{PT1} is precisely an
epimorphism   in $\A^\sharp$ whose image is supported in dimension
$\leq 1$.


\subsection{The stack of framed sheaves}
\label{framed}

 We write $\M$ for the stack of coherent sheaves on $M$. It is an algebraic stack, locally of finite type over $\C$.
  There is another  stack $\M(\OO)$
 with a morphism
 \begin{equation}\label{perv}q\colon \M(\OO)\to \M\end{equation}
 parameterizing coherent sheaves equipped with a section.
 More precisely, the objects of $\M(\OO)$ lying over a scheme $S$ are pairs $(E,\gamma)$ consisting of an $S$-flat coherent sheaf
$E$ on $S\times M$ together with a section
\begin{equation}
\label{bully} \gamma\colon \OO_{S\times M}\to E.\end{equation}
Given a morphism of schemes $f\colon T\to S,$ and an object
$(F,\delta)$ lying over $T$,  a morphism $\theta\colon (F,\delta)\to (E,\gamma)$  lying over
$f$ is an isomorphism
\[\theta\colon f^*(E)\to F\] such that the diagram
\[\begin{CD}f^*(\OO_{S\times M})&@>{f^*(\gamma)}>> &f^*(E)\\ @V{\operatorname{can}}VV && @V\theta VV
\\\OO_{T\times
M} &@>\delta>> & F\end{CD}\] commutes. Here we have taken the usual step of choosing, for each
 morphism of schemes, a pullback of each coherent sheaf on its target. The symbol
$\operatorname{can}$ denotes the canonical isomorphism of
pullbacks.

The stack property for $\M(\OO)$ follows easily from the stack
property for  $\M$. The morphism $q$ is defined by forgetting the
data of the section $\gamma$ in the obvious way.

\begin{lemma}
\label{snow}
The stack $\M(\OO)$ is algebraic and the morphism $q$ is representable and of finite type.
\end{lemma}

\begin{pf}
Consider a scheme $S$ with a morphism $h\colon S\to\M$ corresponding to a family of sheaves $E$ on $S\times M$.
Form the fibre product of stacks
\[\begin{CD}
Z & @>f>> &\M(\OO)\\
 @VpVV & &@VVqV \\
S &@>h>> &\M\end{CD}\] We claim that the stack $Z$ is
representable by a scheme of finite type over $S$. This is enough,
since if we take $S$ to be an atlas for $\M$, the scheme $Z$
becomes an atlas for $\M(\OO)$.

It is easy to see that $Z$ is fibered in sets and corresponds to
a functor
\[Z\colon\big(\sch/S\big)^{\op}\to \Set\]
that sends a morphism $f\colon T\to S$ to the space of sections
\[H^0(T\times M,f^*(E)).\] 
Standard argument allow us to reduce to the case when $S$ is noetherian.
Results of Grothendieck \cite[Theorem 5.8]{Nit} then show that
 this functor is representable by a linear scheme of finite type over $S$.
\end{pf}

 We shall need the
following  result about the fibres of the map $q$.

\begin{lemma}
\label{wet}
There is a stratification of $\M$ by locally-closed substacks \[\M_r\subset \M\]
 such that the objects of $\M_r(\C)$ are coherent sheaves $E\in\A$ with \[\dim_\C
 H^0(M,E)=r.\]
 The pullback of the morphism $q$ to $\M_r$ is a  locally-trivial fibration in the Zariski topology,
with fibre $\C^r$.
\end{lemma}

\begin{pf}
Given a noetherian scheme $S$ and an $S$-flat coherent sheaf $E$ on $S\times
M$, Grothendieck showed \cite[Theorem 5.7]{Nit} that there is a
coherent sheaf $G(E)$ on $S$ together with an isomorphism of
functors
\[\lHom_{\OO_S}(G(E),-)\isom \pi_{S,*}(E\tensor_{\OO_{S \times M}}
 \pi_S^*(-))\colon \Coh(S) \to \Coh(S).\]  The universal
property shows that such a sheaf $G(E)$ is unique up to a unique
isomorphism, and is well-behaved under pullback, in the sense that
if $f\colon T \to S$ is a morphism of schemes, then there is a
natural isomorphism
\[G(f^*(E))\isom f^*(G(E)).\]
It follows from the existence and basic properties of flattening
 stratifications\comment{Needs $S$ noetherian - same issue in next lemma}
  (see e.g. \cite[Theorem 5.13]{Nit}) that there is a stratification of $S$ by locally-closed subschemes $S_r\subset S$ on which $G(E)$ is locally-free of rank r.
 Moreover these stratifications behave well under pullback.
 It follows that the stack $\M$ has a stratification as claimed.

 Using the same notation as in  the proof of Lemma \ref{snow},
 if the morphism $h$ factors via the inclusion $\M_r\subset \M$,
 the stack $Z$ is represented by the total space of the vector bundle
  associated to the locally-free sheaf $G(E)$. It is therefore a
  locally-trivial  fibration with fibre $\C^r$.
\end{pf}


\subsection{Moduli of stable pairs}

We shall view the Hilbert scheme and the moduli space of stable
pairs as open substacks of the moduli stack $\M(\OO)$.

\begin{lemma}
\label{snowy} There are open substacks
\[\Hilb_M\subset\M(\OO),\quad \Hilbs_M\subset\M(\OO),\] whose
$\C$-valued points are morphisms
 \begin{equation}
 \label{nerfititi}
 \gamma\colon \OO_M \to E\end{equation}
 that are  epimorphisms in the category $\A$ or $\A^\sharp$
 respectively.
 \end{lemma}

 \begin{pf}
 We
 just prove the stable pair case; the Hilbert scheme case is easier and well-known. Suppose given a morphism
 \[\gamma\colon \OO_{S\times M} \to E\]
  as in \eqref{bully}. Given a point $s\in S(\C)$ let $\gamma_s\colon \OO_M \to E_s$ denote the pullback of $\gamma$.
  We must show that the set of points $s\in S(\C)$ for which
 \begin{itemize}
 \item[(a)] the sheaf $E_s$ is pure of dimension 1,
 \item[(b)] the sheaf $\coker(\gamma_s)$ is supported in dimension 0,
 \end{itemize}
is the set of $\C$-valued points of an open subscheme of $S$. Standard arguments mean that we can assume $S$ noetherian.

Property (a) is  known to be open (see e.g.
\cite[Prop. 2.3.1]{HL}). For property (b), set
$Q=\coker(\gamma)$. Then
 \[Q_s\isom \coker(\gamma_s)\]
 and so by the  basic properties of flattening
  stratifications (see e.g. \cite[Theorem 5.13]{Nit}) the set of points $s\in S(\C)$ for
   which $\coker(\gamma_s)$ has positive-dimensional support
  is the set of $\C$-valued points of a closed subscheme of $S$.
  \end{pf}

 The substack $\Hilb_M$ is  represented by the Hilbert scheme.
 The substack $\Hilbs_M$ is presumably also represented by a scheme,
 although we have not checked this. Pandharipande and Thomas \cite[Section 1]{PT1} used results of
 Le Potier to show that the open and closed substack
 \[\Hilbs_{M,\leq 1}\subset \Hilbs_M\] parameterizing framed sheaves supported in
 dimension at most 1 is representable.


\section{Behrend function identity}

%

  Recall \cite{Be} that Behrend defined for every scheme $S$ of finite type over $\C$ a constructible function
\[\nu_S\colon S\to\Z,\]
with the property that when $f\colon T\to S$ is a smooth morphism of relative dimension $d$, there is an identity
\begin{equation}\label{snowball}
\nu_T=(-1)^d f^*(\nu_S).\end{equation}
Using this identity it is easy \cite[Prop. 4.4]{JS} to extend Behrend's construction    to give a  locally-constructible
function for any Artin stack locally of finite type over $\C$. The identity \eqref{snowball} then also holds for  smooth morphisms of
stacks.

In this section we establish a Behrend function identity that we
need for the proof of Theorem \ref{conjdt}. The basic idea for the
proof came from unpublished notes of Pandharipande and Thomas;
nonetheless the author takes full responsibility for the details!

\begin{thm}
\label{rahul}
Suppose $\gamma\colon \OO_M\to E$ is a morphism of sheaves such that
\begin{itemize} \item[(a)] $E$ is supported in dimension
$\leq 1$,  \item[(b)]  $\coker(\gamma)$ is supported in dimension
zero.
\end{itemize} Then there is an equality of Behrend functions
\[\nu_{\M(\OO)}(\gamma) = (-1)^{\chi(E)} \cdot \nu_\M(E).\]
\end{thm}

\begin{pf}By Lemma \ref{wet}, the morphism $q$ is a vector bundle over
 the  open and closed substack of $\M$ parameterizing
zero-dimensional sheaves, so if $E$ is supported in dimension 0
the result follows from \eqref{snowball}. Thus we can  assume that
$E$ has positive-dimensional support.

Let $\OO_C\subset E$ be the image of $\gamma$. It is
the structure sheaf of a subscheme $C\subset M$ of dimension $1$.
Take an ample line bundle $L$ such that
\begin{equation}\label{fluss} H^i(M,E\tensor L)=0 \text{ for all }
i>0,\end{equation}
and a divisor $H\in |L|$ such that
 $H$ meets $C$ at finitely many pure-dimension 1 points, not in the
support of $\coker(\gamma)$.

There is a short exact sequence
\[0\lra \OO_M\lRa{s}L\lra \OO_H(H)\lra 0\]
where $s$ is the section of $L$ corresponding to the divisor $H$.
Tensoring it with $E$ and using the above assumptions gives a
 diagram of sheaves
\begin{equation}\label{out}\xymatrix@C=1.5em{  &\OO_{ M} \ar[dr]^{\delta}\ar[d]^{\gamma}
\\
0 \ar[r] &E \ar[r]^{\alpha}& F  \ar[r]^\beta &Q \ar[r] &0 }
\end{equation}
where $F=E\tensor L$.
 The support of the sheaf $Q$ is zero-dimensional,
disjoint from the support of $\coker(\gamma)$ and contained in the
pure dimension 1 part of $C$. In particular,
\begin{equation}
\label{desp} \Hom_M(Q,F)=0.
\end{equation}

Consider  the two points of the stack $\M(\OO)$ corresponding to
the maps \[\gamma\colon \OO_M \to E, \quad \delta\colon \OO_M \to
F.\] The statement of the Theorem holds for the map $\delta$
because Lemma \ref{wet} together with \eqref{fluss} implies that
\[q\colon \M(\OO)\to\M\] is smooth of relative dimension $\chi(F)=H^0(M,F)$
over an open neighbourhood of the point $F\in\M(\C)$. On the other hand, tensoring sheaves with $L$ defines an
automorphism of $\M$, so the Behrend function of $\M$ at the
points corresponding to $E$ and $F$ is equal. Thus to prove the
Theorem it suffices to show that
\[(-1)^{\chi(E)}\cdot \nu_{\M(\OO)}(\gamma) = (-1)^{\chi(F)}
\cdot \nu_{\M(\OO)}(\delta).\]

 Consider the stack $X$ whose
$S$-valued points are diagrams of $S$-flat sheaves on $S\times M$
of the form
\begin{equation}\label{form}\xymatrix@C=1.5em{  &\OO_{S\times M} \ar[dr]^{\delta_S} \ar[d]^{\gamma_S}
\\
0 \ar[r] &E_S \ar[r]^{\alpha_S}& F_S  \ar[r]^{\beta_S} &Q_S \ar[r]
&0 }
\end{equation}
This is an algebraic stack, locally of finite type over $\C$.
The easiest way to see this  is
 to write $X$ as a fibre product in a similar way to the stack $Z$ of
 Lemma \ref{crucial} below.
There are two morphisms
\[p\colon X\to \M(\OO),\quad q\colon X\to\M(\OO),\]
taking such a diagram to the maps $\gamma_S$ and $\delta_S$
respectively. Passing to an open substack of $X$ we can assume
that \eqref{desp} holds at all $\C$-valued points of $X$. It
follows easily that $p$ and $q$ induce injective maps on stabilizer
groups of $\C$-valued points, and hence are representable.

Using \eqref{snowball} it will be enough to show that at the point
 $x\in X(\C)$ corresponding to the diagram \eqref{out}, the
morphisms $p$ and $q$ are smooth of relative dimension $\chi(Q)$
and $0$ respectively. Suppose $f\colon X \to Y$ is a representable
morphism of stacks, with $X$ and $Y$  locally of finite type over
$\C$. Then $f$ is smooth at a point $x\in X(\C)$ precisely if the
following lifting property is satisfied. Let
\[i\colon S \into T\] be a closed embedding of affine schemes,
corresponding to a homomorphism of Artinian local rings
 $A\to A/I$ with $I^2=0$.
Then  given a
 commutative diagram
\[\begin{CD} S & @>s>> & X
\\
@ViVV && @VVfV \\
T & @>t>>& Y\end{CD} \] taking the closed point of $S$ to $x$,
there is a lifting of $t$ to a morphism $u\colon T \to X$
satisfying $f\circ u\isom t$. This follows immediately from the
corresponding characterization of smooth morphisms of schemes
[EGA, IV.17.14.2].

Consider the morphism $q$ first. Suppose  that we are given a
diagram of the form \eqref{form}, whose pullback to the
closed\comment{Why is residue field at closed point $=\C$ - look
in EGA. Also do dotted arrow in diagram} point of $S$ is
isomorphic to the diagram \eqref{out}. Suppose also given a
morphism of $T$-flat sheaves
\[\delta_T\colon \OO_{T\times M} \to F_T\]
on $T\times M$ which when pulled back to $S$  is isomorphic to
$\delta_S$. Then we must include $\delta_T$ in a diagram
\begin{equation}\label{form2}\xymatrix@C=1.5em{  &\OO_{T\times M} \ar[dr]^{\delta_T} \ar[d]^{\gamma_T}
\\
0 \ar[r] &E_T \ar[r]^{\alpha_T}& F_T  \ar[r]^{\beta_T} &Q_T \ar[r]
&0 }
\end{equation}
which becomes isomorphic to \eqref{form} when pulled back to $S$.

Note that it will be enough to prove the existence of such a lift
on an affine open subset $U\subset M$ containing the support of
$Q$. Indeed, we can then take another open subset $V\subset M$
disjoint from the support of $Q$ such that $M=U\cup V$. Over the
 subset $V$ the extension problem is trivial and we can
glue the results to obtain a lift over $M$.

 Let us therefore replace $M$ by $U$. By our assumptions on the
  support of $Q$, we can then assume that the morphism $\gamma$
 is surjective, so that $E\isom \OO_C$, and that the curve $C$ is
of pure dimension 1.  We can further assume that the line bundle
$L$ restricted to $C$ is trivial, and hence that we also have $F\isom \OO_C$.

Let $E_T$ be the image of $\delta_T$ and fill in the diagram
\eqref{form2}. By Lemma \ref{blow} applied with $S=T$ all the
sheaves appearing are flat over $T$.  It is then automatic that
the diagram  becomes  isomorphic to \eqref{form} when pulled back to $S$,
because the morphisms $\beta_S$ and $\beta_T$ are the cokernels of
$\delta_S$ and $\delta_T$ respectively, and the pullback
functor preserves cokernels.
  The lift we
constructed is easily seen to be unique, and hence $q$ is
{\'e}tale.

Let us now check the lifting property for the map $p$. Take a
diagram \eqref{form} again, and suppose this time that we are
given a morphism of $T$-flat sheaves
\[\gamma_T\colon \OO_{T\times M} \to E_T\]
on $T\times M$ which when pulled back to $S$  is isomorphic to
$\gamma_S$. We must include $\gamma_T$ in  a diagram of the form
\eqref{form2} which becomes isomorphic to \eqref{form} when
pulled back to $S$. As before it is enough perform this lifting
over an open affine subset $U\subset M$ containing the support of
$Q$.

By the first part of the proof of Lemma \ref{blow} we can assume
that $E_T \isom \OO_{C_T}$ for some closed subscheme $C_T \subset
T \times U$ which must be flat over $T$. By Lemma \ref{blow} the
sequence in \eqref{form2} takes the form
\[0\lra \OO_{C_S}\lra \OO_{C_S} \lra Q_S \lra 0,\]
 and thus
defines a point of the  relative Hilbert scheme of the
scheme\comment{$C_T/T$ ?} $C_S/S$. This is smooth at the given
point because the obstruction space \cite[Section 5]{Gr} is
\[\Ext^1_C(\OO_C,Q)=0.\]
Thus we can extend the sequence over $T$ to give a sequence
\[0\lra K_T\lra \OO_{C_T} \lra Q_T \lra 0.\]
Applying Lemma \ref{blow} again with $S=T$ shows that $K_T
\isom \OO_{C_T}$ and so we obtain the required extension.
 In the case that \[S=\Spec(\C), \quad T=\Spec
\C[\epsilon]/(\epsilon^2),\] the set of possible extensions is the
tangent space
\[\Hom_C(\OO_C,Q)=\C^{\chi(Q)},\]
and hence $p$ is smooth of this relative dimension.
\end{pf}

\begin{lemma}
\label{blow} Suppose $S$ is the spectrum of an Artinian local
$\C$-algebra with residue field $\C$, and  let $U\subset M$ be an
open affine subscheme.
 Suppose given a short exact sequence of
coherent sheaves
\begin{equation} \label{formy}
 0\lra E_S \lra  F_S\lra Q_S
\lra 0 \end{equation} on $S\times U$, with $F_S$ flat over $S$,
which on restricting to the special fibre gives a sequence of the
form
\[ \cdots\lra \OO_Y\lra  \OO_Z\lra Q \lra 0,\]
with $Z\subset U$ a Cohen-Macaulay curve, and $Q$ supported in
dimension 0. Then all the sheaves appearing in \eqref{formy} are
$S$-flat, and \[E_S \isom F_S \isom \OO_{Z_S}\] for some closed
subscheme $Z_S \subset S \times U$.
\end{lemma}

\begin{pf}
Since $U$ is affine, we can lift the epimorphism $\OO_U\to \OO_Y$
to a morphism $\OO_{S\times U}\to E_S$. This is then surjective by
Nakayama's Lemma. The same argument applies to $F_S$. Thus
\[E_S
\isom \OO_{Y_S}, \quad F_S \isom \OO_{Z_S},\] for subschemes
$Y_S,Z_S \subset S\times U$. The inclusion $E_S \subset F_S$
implies that $Y_S \subset Z_S$ and hence there is a surjection
\begin{equation}
\label{form3} \OO_{Z_S} \onto \OO_{Y_S}.\end{equation} When
restricted to the special fibre this becomes a surjection $\OO_Z
\onto \OO_Y$ which is an isomorphism away from the support of $Q$.
Since $Q$ has dimension zero, and $Z$ is Cohen-Macaulay this must
be an
isomorphism. Thus $Z=Y$. Restricting 
\eqref{formy} to the special fibre, and using purity of $\OO_Y$
 shows that $Q_S$ and hence also $E_S$ are flat over
$S$. If $K_S$ is the kernel of the surjection \eqref{form3}, then
restricting to the special fibre and using flatness of $E_S$ shows
that $K_S=0$.
\end{pf}


\section{Hall algebra identities}
\label{identities} Our strategy for proving Theorem 1.1 is to
first prove  identities in a motivic Hall algebra, and to then
apply an integration map to obtain identities of generating
functions for  DT invariants. In this section we
prove some identities in an infinite-type version of the Hall algebra.

We start by briefly recalling the definition of the motivic Hall
algebra $\RH(\A)$. We refer the reader to \cite[Section 3.4]{B1}
for more details. We then introduce an infinite-type version of
the Hall algebra $\Fi{\A}$. This  provides a simple context in
which to state the identities of this section without having to
worry about convergence issues. Although these identities
themselves will not be used in the rest of the paper, the
arguments used to prove them will appear several times later on.

 One difference from
\cite{B1} which will be important later (see Theorem \ref{deep}),
is that we will always take our Grothendieck groups to be defined
with complex rather than integral coefficients (we could also have
taken rational coefficients).

\subsection{Motivic Hall algebra}
\label{elements}

Suppose $S$ is a stack, locally of finite type over $\C$ and with
affine stabilizers. There is a relative Grothendieck group
$\kst{S}$ defined as the complex vector space spanned by
equivalence classes of symbols of the form
\begin{equation}
\label{symbol} [X\lRa{f}S],\end{equation} with  $X$ a stack of
finite type over $\C$ with affine stabilizers, and $f$ a morphism
of stacks, modulo the following relations \begin{itemize}
\item[(a)] for every pair of stacks $X$ and $Y$ a
relation\[[X_1\amalg X_2\lRa{f_1\sqcup f_2} S]=[X_1\lRa{f_1}
S]+[X_2\lRa{f_2} S],\]

\item[(b)] for every commutative diagram
\[ \xymatrix@C=.8em{
X_1\ar[dr]_{f_1}\ar[rr]^{g} && X_2\ar[dl]^{f_2}\\
&S }
\]with $g$ an equivalence on $\C$-valued points, a  relation
\[[X_1\lRa{f_1} S]=[X_2\lRa{f_2} S],\]

\item[(c)]  for every pair of Zariski locally-trivial fibrations
\[h_1\colon X_1 \to Y, \quad h_2\colon X_2\to Y\]
  with the same fibres, and
every morphism $g\colon Y\to S$, a relation
\[[X_1\lRa{g\circ h_1} S]=[X_2\lRa{g\circ h_2} S].\]
\end{itemize}
The group $\kst{S}$ has the structure of a $\kst{\C}$-module defined by setting
\[[X]\cdot [Y\lRa{f} S]=[X\times Y\lRa{f\circ\pi_2} S]\]
and extending linearly.

The motivic Hall algebra \cite[Section 4]{B1} is the vector space
\[\RH(\A)=\kst{\M}\]
equipped with a non-commutative product  given explicitly  by the
rule
\[ [X_1\lRa{f_1}\M] * [X_2\lRa{f_2} \M] = [Z\lRa{b\circ h}\M], \]
where  $h$ is defined by the following Cartesian square
\begin{equation}
\label{harrogate}\begin{CD}
Z & @>h>> &\M^{(2)} &@>b>> \M\\
 @VVV  &&@VV(a_1,a_2)V \\
X_1\times X_2 &@>f_1\times f_2>> &\M\times\M\end{CD}\end{equation}
It is an algebra over $\kst{\C}$.


\subsection{Infinite-type version}

Given a stack $S$ as before, define an infinite-type Grothendieck
group
 $\operatorname{L(St_{\infty}/S)}$ by considering symbols
\eqref{symbol} as above, but with $X$ only assumed to be locally
of finite type over $\C$. We  also drop the relation (a) since
otherwise  one gets the zero vector space. The infinite-type Hall
algebra is  then the vector space
\[\Fi{\A}=\operatorname{L(St_{\infty}/\A)}\]
  with the convolution product defined exactly as before.

We  shall need some notation for various particular elements of
$\Fi{\A}$. Given an open substack $\N\subset \M$ we write
\[\one{\N} = [\N\lRa{i} \M]\in \Fi{\A},\] where $i\colon \N\to \M$
denotes the inclusion map. Pulling back the morphism
\eqref{perv}
 to  $\N\subset \M$ gives a stack $\N(\OO)$
with a morphism $q\colon \N(\OO)\to \N$, and hence  an element
\[\O{\N}=[\N(\OO)\lRa{q} \M]\in\Fi{\A}.\]

We  abuse notation by using the same symbol for an open substack
of $\M$ and the corresponding full subcategory of $\A$ defined by
its $\C$-valued points. Thus, for example, by \cite[Prop.
2.3.1]{HL}, the objects of the categories $\P$ and $\Q$ of Section
\ref{pairs} are the $\C$-valued points of open substacks
$\P\subset \M$ and $\Q\subset\M$, and there are corresponding
elements \[\one{\P},\one{\Q}\in\Fi{\A}.\]

Finally, restricting $q$
to the open substacks of Lemma \ref{snowy} defines
further elements
\[\H=[\Hilb_M \lRa{q} \M] \in \Fi{\A}, \quad \Hs=[\Hilbs_M \lRa{q} \M] \in \Fi{\A}.\]
We will now establish some identities in the ring $\Fi{\A}$
relating these elements.


\subsection{Torsion pair identities}

\begin{lemma}
\label{torsionpair} There is an identity  \[\one{\M}=\one{\P} *
\one{\Q}.\]
\end{lemma}

\begin{pf}
Form  a Cartesian square
\[\begin{CD}
Z & @>f>> &\M^{(2)}&@>b>>&\M \\
 @VVV & &@VV(a_1,a_2)V \\
\P\times\Q &@>i>> &\M\times\M\end{CD}\] where the map $i$ is the
open embedding. Since the morphism $(a_1,a_2)$ satisfies the
iso-fibration property of \cite[Lemma A.1]{B1}, the groupoid of
$S$-valued points of $Z$ can be described as follows. The objects
are short exact sequences of $S$-flat sheaves on $S\times M$ of
the form
\begin{equation}
\label{bl}0\lra P\lra E\lra Q\lra 0\end{equation} such that $P$
and $Q$ define flat families of sheaves on $M$ lying in the
subcategories $\P$ and $\Q$ respectively. The morphisms are
isomorphisms of short exact sequences. The composition
\[g=b\circ f\colon Z\to \M\] sends such a short exact sequence to
the object $E$. This morphism induces an equivalence\comment{Bit
more explanation} on $\C$-valued points because of the torsion
pair property: every object $E$ of $\A$ fits into a unique short
exact sequence of the form
 \eqref{bl}. Thus the identity follows from relation (b) above.
\end{pf}

\begin{lemma}
\label{harder}
There is an identity \[\O{\M}=\O{\P} * \O{\Q}.\]
\end{lemma}

\begin{pf}
Form Cartesian squares
\[\begin{CD}
Y  @>p>> &X  @>j>> &\M^{(2)}&@>b>>&\M\\
@VVV & @VVV & @VV(a_1,a_2)V \\
\P(\OO) \times \Q(\OO)@>(q,\id)>> &\P\times\Q(\OO) @>(i,q)>>
&\M\times\M\end{CD}\] Then $\O{\P}*\O{Q}$ is represented by the
composite map $b\circ j\circ p\colon Y\to \M$. Note that, by Lemma \ref{wet},
the map
\[q\colon \P(\OO)\to\P\]
is a Zariski locally-trivial  fibration, with fibre over a zero-dimensional sheaf $P$ being the vector space $H^0(M,P)$.
By pullback the same is true of the morphism $p$.

 The groupoid of $S$-valued points
of $X$ is as follows. The objects are short exact sequences of
$S$-flat sheaves on $S\times M$ of the form
\[0\lra P\lRa{\alpha} E\lRa{\beta} Q\lra 0\]
such that $P$ and $Q$ define flat families of sheaves on $M$ lying
in the subcategories $\P$ and $\Q$ respectively, together with a
map \[\gamma\colon \OO_{S\times M} \to Q.\] The morphisms are
isomorphisms of short exact sequences commuting with the map $\gamma$.
We can represent these objects as diagrams
\begin{equation}\label{bloo}\xymatrix@C=1.5em{  &&&\OO_{S\times M} \ar[d]^{\gamma}
\\
0 \ar[r] &P \ar[r]_{\alpha}& E  \ar[r]_\beta &Q \ar[r] &0 }
\end{equation}

Consider the stack $Z$ of Lemma \ref{torsionpair} with its map $g\colon Z
\to\M$ and form a Cartesian square
\[\begin{CD}
W & @>h>> &\M(\OO)\\
 @VVV & &@Vq VV \\
Z &@>g>> &\M\end{CD}\]
 Since $g$ induces an equivalence on $\C$-valued points, so
too does $h$, so  the element $\O{\M}$ can be represented by
the morphism $q\circ h$.

The groupoid of $S$-valued points of $W$ can be represented by diagrams
\begin{equation}
\label{needle}\xymatrix@C=1.5em{  &&\OO_{S\times M} \ar[d]^\delta
\\
0 \ar[r] &P \ar[r]_{\alpha}& E  \ar[r]_\beta &Q \ar[r] &0 }
\end{equation}
Setting $\gamma=\beta\circ \delta$ defines a morphism of stacks
$W\to X$. It is easy to see that this is a Zariski
locally-trivial fibration, with fibre over a $\C$-valued point of
$X$ represented by a diagram \eqref{bloo} being an affine space
for $H^0(M,P)$. Since this is the same fibre as the map $p$ the
result follows from relation (c) above.
\end{pf}


\subsection{Hilbert scheme identity}

The next result is analogous to an identity proved by Engel and Reineke \cite[Lemma 5.1]{ER} in the context of Hall algebras of quiver representations
defined over finite fields.

\begin{lemma}
\label{crucial} There is an identity \[\O{\M} = \H
* \one{\M}.\]
\end{lemma}

\begin{pf}
Form a Cartesian square
\[\begin{CD}
Z & @>s>> &\M^{(2)}&@>b>>&\M\\
 @V(r,a_2\circ s)VV & &@VV(a_1,a_2)V \\
\Hilb_M\times\M &@>(q,\id)>> &\M\times\M\end{CD}\] 
The groupoid of $S$-valued points of $Z$ are as follows. The objects are
short exact sequences
of $S$-flat sheaves on $S\times M$
\[0\lra A\lRa{\alpha} E\lRa{\beta} B\lra 0\]
together with an epimorphism $\gamma\colon \OO_{S\times M} \onto
A$. The morphisms are isomorphisms of such short exact sequences
commuting with the maps  from $\OO_{S\times M}$.
We can represent these objects as diagrams of the form
\begin{equation}\label{diaia}\xymatrix@C=1.5em{  &\OO_{S\times M} \ar[d]^{\gamma}
\\
0 \ar[r] &A \ar[r]^{\alpha}& E  \ar[r]^\beta &B \ar[r] &0 }
\end{equation}

There is a morphism of stacks $h\colon Z\to \M(\OO)$ sending a diagram \eqref{diaia} to the composite map
 \[\delta=\alpha\circ\gamma\colon\OO_{S\times M}\to E.\]
  This morphism fits
 into
a commuting diagram of stacks
\begin{equation} \xymatrix@C=1em{
Z\ar[rr]^{h} \ar[dr]_{b\circ s}&& \M(\OO)\ar[dl]^{q}\\  &\M }
\end{equation}
Then $h$ induces an equivalence on $\C$-valued points\comment{Bit
more explanation}, because every map $\delta\colon \OO_M\to E$ in
$\A$ factors uniquely as an epimorphism $\gamma\colon \OO_M\onto
A$ followed by a monomorphism $\alpha\colon A\into E$.
\end{pf}


\subsection{Stable pair identity}


\begin{lemma}
\label{stablepairs}There is an identity
\[\O{\Q}=\Hs*\one{\Q}.\]
\end{lemma}

\begin{pf}
This is similar to Lemma \ref{crucial}. We consider the Cartesian
square
\[\begin{CD}
Z & @>s>> &\M^{(2)}&@>b>>&\M\\
 @V(r,a_2\circ s)VV & &@VV(a_1,a_2)V \\
\Hilbs_M\times\Q &@>(q,i)>> &\M\times\M\end{CD}\] This time the
objects in the groupoid of $S$-valued points of $W$ are diagrams
of $S$-flat sheaves on $S\times M$ of the form
\begin{equation}\xymatrix@C=1.5em{  &\OO_{S\times M} \ar[d]^{\gamma}
\\
0 \ar[r] &A \ar[r]^{\alpha}& E  \ar[r]^\beta &B \ar[r] &0 }
\end{equation}
such that for each point $s\in S(\C)$, the  pullback $B_s$ is in
$\Q$, and the map $\gamma_s\colon \OO_M\to A_s$
 is an
epimorphism in $\A^\sharp$. In particular, by Lemma \ref{oneup} this implies that $A_s$ is
in $\Q$ and hence so too is $E_s$.

There is a morphism $h\colon Z\to \Q(\OO)$ sending such an object
to the composite morphism $\delta=\alpha\circ\gamma$ as before,
and a commutative diagram
\begin{equation} \xymatrix@C=1em{
Z\ar[rr]^{h} \ar[dr]_{b\circ s}&& \Q(\OO)\ar[dl]^{q}\\  &\M }
\end{equation}
Once again $h$ is an equivalence on $\C$-valued points, this time because  if $E\in\Q$
then every map $\OO_M\to E$ factors uniquely as an epimorphism
$\OO_M \to A$ in $\A^\sharp$ followed by a monomorphism $A\to E$
in $\A^\sharp$. By Lemma \ref{oneup}, one then  has $A\in\Q$, and forming the short exact sequence
\[0\lra A\lra E\lra B\lra 0\]
in $\A^\sharp$, and applying the argument of Lemma \ref{oneup}
gives $B\in\Q$ also.
\end{pf}


\section{Hall algebras of dimension one sheaves}

In this section we assemble some general properties of Hall algebras 
of coherent sheaves on $M$ supported in dimension $\leq 1$. We first adapt some material from \cite{B1} to this context, and then introduce a formal completion operation which
allows us to treat the infinite sums of invariants appearing in Theorem \ref{conjdt}. Finally, we consider the properties of the dualizing functor applied to pure dimension 1 sheaves.

\subsection{Hall algebra and integration map}
We write
\[\CC=\Coh_{\leq 1}(M)\subset\Coh(M)=\A\]
for the full subcategory consisting of sheaves
 with support of dimension $\leq 1$.  By our usual abuse of
notation, this corresponds to an open and closed substack
\[\CC\subset \M.\] As in \cite[Section 4.4]{B1} there is a
$\kst{\CC}$-subalgebra $\RH(\CC)\subset \RH(\A)$ spanned by symbols
\begin{equation} \label{sym} [X\lRa{f}\M]\end{equation} factoring
via the inclusion  $\CC\subset \M$. As explained in \cite[Section
5.1]{B1} there is
 a $\kvar{\C}[\LL^{-1}]$-subalgebra
\[\RHreg(\CC)\subset \RH(\CC)\]
 spanned by  symbols \eqref{sym} with $X$ a scheme. Elements of this subalgebra are called regular.
The  quotient algebra
\[\RHsc(\CC)=\RHreg(\CC)/(\LL-1)\RHreg(\CC)\]
is commutative, and is equipped with a Poisson bracket
\[\{f,g\}=\frac{f*g-g*f}{\LL-1} \mod (\LL-1).\]
All these algebras are  graded by the commutative monoid $\Delta$,
the basic reason being the decomposition
\begin{equation}
\label{yehbaby!}\CC=\big.\coprod_{\alpha\in \Delta}
\CC_\alpha\end{equation} into open and closed substacks
paramaterizing sheaves with a fixed Chern character.

 Let
$\C[\Delta]$ denote the monoid algebra of $\Delta$ over $\C$. It
has basis given by symbols $x^\alpha$ for $\alpha\in\Delta$ with
product defined by
\[x^\alpha*x^\beta=x^{\alpha+\beta}.\]
Since the Euler form on $N_{\leq 1}(M)$ is trivial we equip $\C[\Delta]$
 with the trivial Poisson bracket.
 \begin{thm}
 There is a
homomorphism of $\Delta$-graded Poisson algebras
\[I\colon \RHsc(\CC)\to \C[\Delta]\]
with the property that
\begin{equation}
\label{rich} I\big([X\lRa{f}\CC_\alpha]\big)=\chi(X,f^*(\nu_{\CC}))
\cdot x^\alpha,\end{equation} where $\nu_{\CC}\colon \CC\to\Z$ is the
Behrend function for the  stack $\CC$.
\end{thm}

\begin{pf}
This follows immediately from \cite[Theorem 5.2]{B1} upon
restricting to the subalgebra $\RH(\CC)\subset \RH(\A)$.
\end{pf}


\subsection{Laurent subsets}
 Theorem \ref{conjdt} involves Laurent series rather than polynomials. To fit these into our framework we need to
 enlarge the algebras $\RH(\CC)$ and $\C[\Delta]$ in a formal way using the grading.
 In this section we explain how to do this.

\begin{defn}A subset $S\subset \Delta$  will be called Laurent if
for each $\beta\in N_1(M)$ the set of integers $n$ for which
$(\beta,n)\in S$ is bounded below.
\end{defn}

Let $\Phi$ denote the set of Laurent subsets of $\Delta$. By Lemma \ref{eff} this
system of subsets satisfies
\begin{itemize}
\item[(a)] if $S,T \in\Phi$ then so is $S+T=\{\alpha+\beta:
\alpha\in S, \beta\in T\}$.

\item[(b)] if $S,T\in\Phi$ and $\alpha\in \Delta$ there are only
finitely many decompositions $\alpha=\beta+\gamma$ with $\beta\in
S$ and $\gamma\in T$.
\end{itemize}

 Suppose that $A=\bigoplus_{\gamma\in\Delta} A_\gamma$ is a $\Delta$-graded ring.  We can form an abelian group  $A_\Phi$ whose
elements are formal sums
\begin{equation}
\label{peter}
a=\sum_{\gamma\in S\subset \Delta} a_\gamma\end{equation}
for Laurent subsets $S\in\Phi$. Given  an element $a\in A_\Phi$ as in \eqref{peter}, we
set $\pi_\gamma({a})=a_\gamma\in A_\gamma$. The ring structure on
$A$ induces one on $A_\Phi$ via the rule
 \[\pi_\gamma(a*b)=\sum_{\alpha+\beta=\gamma} \pi_\alpha(a) *
 \pi_\beta(b).\]
We can equip $A_\Phi$ with a topology by stipulating that
a sequence $(a_j)\subset A_\Phi$ is convergent if for any
$(\beta,m)\in \Delta$ there is an integer $K$  such that
 \[i,j\geq K \implies \pi_{(\beta,n)}(a_i)= \pi_{(\beta,n)}(a_j)\text{ for all }n \leq m. \]
 It is easy to check using Lemma \ref{eff} that with this topology $A_\Phi$ becomes a
topological algebra.

\begin{lemma}
\label{invert}
If $A$ is a $\Delta$-graded $\C$-algebra, and $a\in A_\Phi$ satisfies $\pi_0(a)=0$,
then any series
\[\sum_{j\geq 1} c_j a^j,\]
with coefficients $c_j\in\C$ is convergent in the topological ring
$A_\Phi$. In particular the element $1-a$ is invertible.
\end{lemma}

\begin{pf}Define a subset
\[\Delta(a)=\{\gamma\in\Delta : \pi_\gamma(a)\neq 0\}.\]
Then $\Delta(a)$ is a Laurent subset consisting of nonzero elements. Fix $(\beta,m)\in\Delta$ and consider elements $\gamma=(\beta,n)\in\Delta$ with $n\leq m$.
It follows from Lemma \ref{eff} that there are only finitely many decompositions of such an element $\gamma$ into a sum of the form
\[\gamma=\gamma_1+\cdots +\gamma_k, \quad \gamma_i\in \Delta(a).\]
The convergence of the given series  follows immediately from this.
 The final statement follows by considering the geometric series
in the usual way.
\end{pf}

If $A$ and $B$
are two $\Delta$-graded algebras then a $\Delta$-graded
homomorphism $f \colon A\to B$ induces a
 continuous homomorphism $f_\Phi\colon A_\Phi\to B_\Phi$ by
setting
\[\pi_\gamma(f_\Phi(a))= f(\pi_\gamma(a)).\]
Note that if $f$ is injective then so is $f_\Phi$. Thus if $A$ is
a subring of $B$ we can identify $A_\Phi$ with a subring of
$B_\Phi$.


\subsection{$\Phi$-finite elements}

Applying the construction of the last section to the Hall algebra
gives a $\kst{\C}$-algebra $\RH(\CC)_\Phi$ with a
$\kvar{\C}[\LL^{-1}]$-subalgebra
\[\RHreg(\CC)_\Phi\subset \RH(\CC)_\Phi\]
whose elements we again call regular. There is also a commutative
quotient algebra
\[\RHsc(\CC)_\Phi=\RHreg(\CC)_\Phi/(\LL-1)\RHreg(\CC)_\Phi\]
equipped with a Poisson bracket, and
 a  homomorphism of Poisson algebras
\begin{equation}
\label{harrow}I_\Phi\colon \RHsc(\CC)_\Phi \to
\C[\Delta]_{\Phi}.\end{equation}

Recall the decomposition \eqref{yehbaby!} of the stack $\CC$.


\begin{defn}
A morphism of stacks  $f\colon X \to \CC$  will be called
$\Phi$-finite if
\begin{itemize}
\item[(a)] $X_\alpha=f^{-1}(\CC_\alpha)$ is of finite type for all
$\alpha\in\Delta$, 
\item[(b)] there is a  subset $S\in\Phi$ such that $X_\alpha$ is
empty unless $\alpha\in S$.  \end{itemize}
\end{defn}

Any such morphism defines an element of the ring $\RH(\CC)_\Phi$
defined by the formal sum
\[\sum_{\alpha\in S}[X_\alpha\lRa{f} \CC].\]
To give some examples,  introduce the open and closed subschemes
\[\Hilb_{M,\leq 1}=\CC \cap \Hilb_M, \quad \Hilbs_{M,\leq 1}=\CC \cap
\Hilbs_M\] parameterizing quotients of $\OO_M$ supported in
dimension $\leq 1$.

\begin{lemma}
\label{john}
 The  morphisms
\[q\colon \Hilb_{M,\leq 1}\to \CC, \quad q\colon \Hilbs_{M,\leq
1}\to \CC,\] are $\Phi$-finite. The corresponding elements
$\H_{\leq 1}$ and $\Hs_{\leq 1}$ of $\RH(\CC)_\Phi$ satisfy
\begin{equation}
\label{john1} I_\Phi(\H_{\leq 1})=\sum_{(\beta,n)\in \Delta}
(-1)^n \DT(\beta,n) x^\beta q^n=\DT_\beta(-q)\cdot x^\beta,
\end{equation} where we have written $x^\beta=x^{(\beta,0)}$ and $q=x^{(0,1)}$.
Similarly
\begin{equation}
\label{john2} I_\Phi(\Hs_{\leq 1})=\sum_{(\beta,n)\in \Delta}
(-1)^n \PT(\beta,n) x^\beta q^n=\PT_\beta(-q)\cdot x^\beta.\end{equation}
\end{lemma}

\begin{pf}The Hilbert scheme parameterising quotients of
$\OO_M$ of fixed class $\gamma\in N(M)$ is always of finite type.
On the other hand the set of elements $\gamma\in \Delta$ for which
$\Hilb_M(\gamma)$ is non-empty is Laurent, the basic reason being
that the curves lying in a fixed class $\beta\in N_1(M)$ have
bounded genus.   The same argument applies to the stable pair
moduli space.  The formulae then follow from  Theorem \ref{rahul}
and Behrend's description of DT
invariants as a weighted Euler characteristic.
\end{pf}




\subsection{Duality functor}
\label{dualt}

There is a full subcategory $\Qlast\subset \CC$
consisting of pure dimension 1 sheaves.

\begin{lemma}
The functor
\[\DD=\lExt^2_{\OO_{M}}(-,\OO_M)\colon \Coh(M) \to \Coh(M)\]
restricts to a contravariant equivalence
\[\DD\colon\Qlast\to \Qlast\]
satisfying $\DD^2\isom \id$.
\end{lemma}

\begin{pf}
It is immediate from the definition that the functor
\[\DD=\R\lHom_{\OO_M}(-,\OO_M)[2]\colon \D^b\Coh(M) \to \D^b\Coh(M)\] is a contravariant  equivalence,
and satisfies $\DD^2\isom \id$.
Thus we must just show that with the given shift $\DD$ takes the
subcategory $\Qlast$ into itself.

By Serre duality, for any object $E\in \Qlast$ and any point $x\in
M$
\begin{equation}
\label{burble} \Ext^3_M(E,\OO_x)=\Hom_M(\OO_x,E)^* =
0,\end{equation} and so by the usual argument (see e.g.
\cite[Prop. 5.4]{BM}) $E$ has a locally-free resolution of length
$2$. Applying $\DD$  gives another length $2$ complex of
locally-free sheaves $\DD(E)$ with the same support as $E$. By a
standard  argument (see e.g. \cite[Lemma 4.2]{BM}) it follows that
$\DD(E)$ is concentrated in degree $0$, and \eqref{burble} then
shows that $\DD(E)$ is pure.
\end{pf}

The functor $\DD$ preserves families and hence induces an
involution $\DD$ of the open substack $\Qlast\subset \CC$. There
is a $\kst{\C}$-submodule   \[\kst{\Qlast} \subset \kst{\CC}\]
spanned by maps \eqref{symbol} factoring via the open substack
$\Qlast$. The involution $\DD$  defines  a $\kst{\C}$-linear
involution
\[\DD_* \colon \kst{\Qlast}\to \kst{\Qlast}.\] This will ultimately be
responsible for the invariance of the DT generating function under
the transformation $q\leftrightarrow q^{-1}$.

\begin{lemma}
The submodule   $\kst{\Qlast}$ is closed under the convolution
product. It therefore defines a subalgebra $\RH(\Qlast) \subset
\RH(\CC)$. Moreover,  one has\[\DD_*(a
* b)=\DD_*(b)
* \DD_*(a),\] for any $a,b\in \RH(\Qlast)$, or in other words, $\DD_*$ is an anti-involution of
$\RH(\Qlast)$.
\end{lemma}

\begin{pf}
The first statement is just the fact that the subcategory $\Qlast$ is
closed under extensions. The second statement holds because $\DD$
is an exact anti-equivalence, which means that if
\[0\lra A_1\lra B\lra A_2\lra 0\]
is a short exact sequence in $\Qlast$, then so is
\[0\lra \DD(A_2)\lra \DD(B) \lra \DD(A_1)\lra 0.\]
Thus $\DD$ induces an automorphism of the stack of short exact
sequences in $\Qlast$ permuting the two projections $a_1,a_2$ to $\Qlast$.
\end{pf}


\section{DT/PT correspondence}

In this section we introduce a  stability condition on the
category $\CC$ and state a deep result of Joyce
which shows that certain elements relating to characteristic
functions of semistable objects are regular.  We also
prove a Hall algebra identity arising from the existence and
uniqueness of Harder-Narasimhan filtrations. These results are then used to prove the first part of Theorem  \ref{conjdt}, the correspondence between Donalson-Thomas
and stable pair invariants.

\subsection{A stability condition}
\label{ex}

Choose an ample divisor $H$ on $M$. Given a class $\gamma\in \Delta$ define
the slope
\[\mu(\gamma)=\frac{\ch_3(\gamma)}{\ch_{2}(\gamma)\cdot H}\in (-\infty,+\infty].\]
If $\ch_{2}(\gamma)=0$ we consider $\gamma$ to have slope
$+\infty$, otherwise $\mu(\gamma)\in\QQ$. A nonzero object $E\in
\CC$ is (Gieseker or Simpson) semistable  if \[\mu(A)\leq \mu(E)\]
for all nonzero subobjects $A\subset E$.  We write
$\semi(\gamma)\subset \CC$ for the open stack whose $\C$-valued
points are semistable sheaves of class $\gamma\in \Delta$.

Every nonzero sheaf $E\in\CC$ has a unique Harder-Narasimhan
filtration
\[0=E_0\subset E_1\subset \cdots\subset E_{n-1}\subset E_n=E,\]
whose factors $F_i=E_i/E_{i-1}$ are semistable with descending slope
\[\mu(F_1)> \mu(F_2)>\cdots >\mu(F_n).\]
Given an interval $I\subset (-\infty,+\infty]$ define
$\semi(I)\subset \CC$ be the full subcategory consisting of zero
objects together with those sheaves whose Harder-Narasimhan
factors all have slope in $I$. Note that there are identifications
\[\P=\semi(\infty), \quad \Qlast=\semi(-\infty,+\infty).\]

By the argument of \cite[Prop. 2.3.1]{HL}, bounding slopes of
Harder-Narasimhan factors is an open condition, so for each
interval $I$ there is an open substack $\semi(I) \subset \M$. If
$I$ is bounded below, then by \cite[Theorem 3.3.7]{HL} the
inclusion functor is $\Phi$-finite, and so defines an element
\[\one{\semi(I)}\in {\RH(\CC)_\Phi}.\]
In particular, for each  $\mu\in (-\infty,+\infty]$ there is a
full subcategory $\semi(\mu)\subset \CC$ consisting of the zero
objects and the semistable objects of slope $\mu$.
 Since the identity
in the Hall algebra is represented by the substack $[\M_0\subset
\M]$ of zero objects, one has \begin{equation} \label{peasy}
\one{\semi(\mu)} =1 + \sum_{\stackrel{0\neq \gamma\in
\Delta}{\mu(\gamma)=\mu}} \one{\semi(\gamma)}.\end{equation}

 It is easy to see
that for any sheaf $E\in \Qlast$
\begin{equation}
\label{burb} \ch(\DD(E))=
(\ch_2(E),-\ch_3(E))\in\Delta.\end{equation} which of course
implies that
\begin{equation}
\label{burb2} \mu(\DD(E))= -\mu(E).\end{equation} The following is
an easy consequence of this.

\begin{lemma}
\label{new} Let $I\subset \R$ be a bounded interval. Then
\[\DD(\semi(I))=\semi(-I).\]
\end{lemma}

\begin{pf}
Suppose $E\in \Qlast$. Then $E$ lies in the subcategory
$E\in\semi([a,b])$ precisely if for every short exact sequence of
objects of $\Qlast$
\begin{equation}
\label{ben} 0\lra A\lra E\lra B\lra 0\end{equation}
 one has $\mu(A)\leq b$ and $\mu(B)\geq a$.
The result therefore follows from \eqref{burb2}.
\end{pf}


\subsection{Harder-Narasimhan identity}
The next result is analogous to an identity proved by  Reineke \cite[Prop. 4.8]{rei} in the context of Hall algebras of quiver representations defined over finite fields.
Similar identities play a fundamental role in the study of wall-crossing behaviour in  \cite{Jo4,KS}.

\begin{lemma}
\label{HN} Assume the interval $I\subset (-\infty,+\infty]$ is bounded below. Then
 there is an identity
\[\one{\semi(I)}=\prod_{\mu\in I} \one{\semi(\mu)}\in \RH(\CC)_\Phi,\]
 where the infinite product is taken in
descending order of slope.
\end{lemma}

\begin{pf}
We first explain what the identity means.  Given a finite subset
$V\subset I$ we can form a product
\[\one{\semi(V)}=\prod_{\mu\in V} \ \one{\semi(\mu)}\in\RH(\CC)_\Phi,\]
where we take the terms in descending order of slope. Suppose
that
\[V_1\subset V_2\subset \cdots \subset V_j\subset \cdots \subset
I\] is an increasing sequence of subsets whose union contains all
rational points of $I$. Then the claim is that the elements
$\one{\semi(V_j)}$ converge to $\one{\semi(I)}$.

Fix $(\beta,m)\in \Delta$ and consider elements $\gamma=(\beta,n)$
with $n<m$. By the boundedness assumption on the interval $I$, and
Lemma \ref{eff}, there are only finitely many ways of writing such
an element as a sum \begin{equation}
\label{reallybored}\gamma=\gamma_1+ \cdots \gamma_k,\end{equation}
with each $\gamma_i\in \Delta$ and satisfying $\mu(\gamma_i)\in
I$.

 The existence and uniqueness of the Harder-Narasimhan
filtration together
  with the formula for repeated products in the Hall algebra (see \cite[Lemma 4.2]{B1})
   implies that  there is an identity
 \[\pi_\gamma(\one{\semi(I)}) = \sum_{k\geq 1}
  \sum_{\stackrel{\gamma_1+\cdots +\gamma_k=\gamma}{\stackrel{\mu(\gamma_1)>\cdots>\mu(\gamma_k)}{\mu(\gamma_i)\in
  I}}}
  \one{\semi(\gamma_1)} *\cdots *\one{\semi(\gamma_k)},\]
  where the sum on the right is finite. Indeed, each term on the
  right is represented by a stack whose $\C$-valued points
  paramaterize sheaves $E$ with a Harder-Narasimhan filtration of
  the given type.
  By \eqref{peasy} it follows that
  \[\pi_\gamma(\one{\semi(I)})=\pi_\gamma\big(\prod_{\mu\in I} \
\one{\semi(\mu)}\big).\]
 Thus if $j$ is large enough that $V_j$ contains the slopes of all $\gamma_i$ appearing in decompositions of the form
\eqref{reallybored} then
\[\pi_\gamma(\one{\semi(I)})=\pi_\gamma(\one{\semi(V_j)})\] which proves the claim.
\end{pf}


%
\subsection{The no-pole theorem}

 The following
statement is a consequence of deep results of Joyce. It is an
analogue of the no-poles conjecture in \cite{KS}. At this point it becomes
important that we have used complex (or at least rational)
coefficients to define our Grothendieck rings (rather than the
integral coefficients of \cite{B1}).

\begin{thm}
\label{deep}
For each $\mu\in (-\infty,+\infty]$ we can write\footnote{This result is false as stated. It becomes true after further localisation of the ring $\RHreg(\CC)_\Phi$, and the rest of the argument is unaffected. See Section \ref{for} below.}
\[\one{\semi(\mu)}=\exp(\epsilon_\mu)\in {\RH(\CC)_\Phi}\]
with $\eta_\mu=[\C^*] \cdot \epsilon_\mu\in \RHreg(\CC)_\Phi$ a regular element.
\end{thm}

%
\begin{pf}
Given $\gamma\in \Delta$, consider the finite sum
\[\epsilon_\gamma=\sum_{k\geq 1} \sum_{\stackrel{\gamma=\gamma_1+\cdots + \gamma_k}{
\mu(\gamma_i)=\mu(\gamma)}} \frac{(-1)^{k-1}}{k}
\one{\semi(\gamma_1)}
* \cdots
* \one{\semi(\gamma_k)}.\]
Joyce proved  that  $[\C^*]\cdot \epsilon_\gamma\in \RHreg(\CC)$.
By \eqref{peasy} this can be rephrased as the statement that
\[[\C^*]\cdot \log (\one{\semi(\mu)})\in  \RHreg(\CC)_\Phi,\]
which is equivalent to the claim made in the statement of the
Theorem via the inverse properties of the expansions of log and
exp.

More precisely \cite[Theorem 8.7]{Jo3} shows that
$\epsilon_\gamma\in\Pi_1\SF(\M)$, where $\Pi_1$ is Joyce's
projection onto virtual indecomposables. But using
 \cite[Cor. 5.10]{Jo2} and the
definition of the  operator $\Pi_1$ this means that in the stack
functions algebra $\bar{\SF}(\M,\Upsilon,\Lambda)$ the element
$\epsilon_\gamma$ is represented by a sum of elements of the form
\[[X_i/\C^* \lRa{f_i} \M]\]
with each $X_i$ a variety. Taking the motivic invariant to be the canonical map
\[\Upsilon\colon \kvar{\C} \to \kst{\C},\]
and using \cite[Remark 3.11]{B1}, shows that $[\C^*]\cdot
\epsilon_\gamma$ is regular as claimed.
\end{pf}

Joyce's proof of this crucial result  uses the full force of the
technology developed in the long series of papers
\cite{Jo0,Jo1,Jo2}, including moduli spaces of configurations and
the virtual projection operators $\Pi_k$. In the author's opinion
it would be a worthwhile project to try to write down a more
conceptual and easily-understood proof.

\begin{cor}
\label{ad} For any $\mu\in (-\infty,+\infty]$  the element
$\one{\semi(\mu)}\in \RH(\CC)_\Phi$ is invertible, and the
automorphism
\[\add_{\one{\semi(\mu)}} \colon {\RH(\CC)_\Phi} \to {\RH(\CC)_\Phi}\]
preserves the subring of regular elements. The induced Poisson
automorphism of $\RHsc(\CC)_\Phi$ is given by
\[\add_{\one{\semi(\mu)}} = \exp \{\eta_\mu,-\}.\]
\end{cor}

\begin{pf}
That $\one{\semi(\mu)}$ is invertible follows from Lemma
\ref{invert}. The identity\comment{Give proof}
\[\add_{\exp(x)} = \exp (\ad_x)\]
 shows that
\[\add_{\one{\semi(\mu)}}=\exp (\ad_{\epsilon_\mu})=\exp  (\ad_{(\LL-1)^{-1} \eta_\mu}).\]
Theorem \ref{deep} shows that $\eta_\mu$ is regular, and since, by
\cite[Theorem 5.1]{B1} the multiplication on $\RHreg(\A)_\Phi$ is
commutative modulo $(\LL-1)$, it follows that this operation
preserves $\RHreg(\A)_\Phi$. The last statement follows from the
definition of the Poisson bracket on $\RHsc(\A)$.
\end{pf}


\subsection{The DT/PT correspondence}

 Consider the part
of the Hilbert scheme $\Hilb_{M,0}$ parameterizing
zero-dimensional subschemes of $M$. The morphism \[q\colon
\Hilb_{M,0} \to \CC\] is $\Phi$-finite, and so defines an
element $\H_0\in \RH(\CC)_\Phi$. As in Lemma \ref{john} one has
\begin{equation}
\label{newdtwo}I_\Phi(\H_{0})=\sum_{(\beta,n)\in \Delta} (-1)^n
\DT(0,n)  q^n=\DT_0(-q),\end{equation} where  $q=x^{(0,1)}\in\C[\Delta]$.

\begin{prop}
\label{toi}
There is an identity
\[\H_{\leq 1} * \one{\P} =
  \H_0 * \one{\P} * \Hs_{\leq 1} \]
  in ${\RH(\CC)_\Phi}$.
  \end{prop}

  \begin{pf}
We work in the ring ${\RH(\CC)_\Phi}$ throughout. The first claim
is that
  \begin{equation}
  \label{left}\H_{\leq 1} * \one{\semi{([\mu,\infty])}} - \O{\semi{([\mu,\infty])}} \to 0 \text{ as }\mu\to-\infty.\end{equation}
  To see this
 fix $(\beta,m)\in \Delta$ and consider classes $(\beta,n)$ with $n<m$. There are only finitely many decompositions
 of such elements \[(\beta,n)=(\beta_1,n_1) + (\beta_2,n_2)\]
 such that the elements
 \[\pi_{(\beta_1,n_1)}(\H_{\leq 1}), \quad \pi_{(\beta_2,n_2)} (\one{\semi{([\mu,\infty])}}),\]
 are nonzero. By boundedness of
 the Hilbert scheme we can assume that $\mu$ is small enough so that for any of these,  all points $\OO_M\onto A$ of $\Hilb_M(\beta_1,n_1)$ satisfy
 \[A\in \semi{([\mu,\infty])}.\]
Suppose given a diagram of sheaves on $M$ of the form
 \begin{equation}
 \xymatrix@C=1.5em{  &\OO_{M} \ar[d]^{\gamma}
\\
0 \ar[r] &A \ar[r]^{\alpha}& E  \ar[r]^\beta &B \ar[r] &0 }
\end{equation}
with $\gamma$ surjective, and $[E]=(\beta,n)$. Then the above assumption on $\mu$ implies that
\[B\in\semi{([\mu,\infty])} \iff E\in \semi{([\mu,\infty])}.\]
 The argument of Lemma \ref{crucial} then shows that
 \[\pi_{(\beta,n)}\big(\H_{\leq 1} * \one{\semi{([\mu,\infty])}}\big) = \pi_{(\beta,n)}\big(\O{\semi{([\mu,\infty])}}\big)\]
 which proves the claim.
The same argument, but using Lemma \ref{stablepairs} gives
\begin{equation}
\label{right} \Hs_{\leq 1} *
\one{\semi{([\mu,\infty))}}-\O{\semi{([\mu,\infty))}}\to 0\text{
as }\mu\to -\infty.\end{equation} Since $\P=\semi(\infty)$ the
proofs of Lemma \ref{torsionpair} and Lemma \ref{harder} easily
extend to give identities
\[\one{\semi{([\mu,\infty])}}=\one{\P} * \one{\semi{([\mu,\infty))}}, \quad \O{\semi{([\mu,\infty])}}=\O{\P} * \O{\semi{([\mu,\infty))}}.\]
Expression  \eqref{left}  can thus be rewritten
  \[\H_{\leq 1} * \one{\P} * \one{\semi{([\mu,\infty))}}- \O{\P}*\O{\semi{([\mu,\infty))}} \to 0 \text{ as }\mu\to -\infty.\]
  The argument of Lemma \ref{crucial} gives an identity
\begin{equation}
\label{nearlydone}
\O{\P} = \H_0 * \one{\P}.
\end{equation}
  Multiplying  \eqref{right} on the left by $\O{\P}$ and using \eqref{nearlydone} gives
  \[\H_0 * \one{\P} * \Hs_{\leq 1} * \one{\semi{([\mu,\infty))}}-\O{\P}*\O{\semi{([\mu,\infty))}}\to 0 \text{ as }\mu\to -\infty.\]
  Thus
  \[\H_{\leq 1} * \one{\P} * \one{\semi{([\mu,\infty))}}-
  \H_0 * \one{\P} * \Hs_{\leq 1} * \one{\semi{([\mu,\infty))}}\to 0 \text{ as }\mu\to -\infty.\]
 By Lemma \ref{invert} the element $\one{\semi{([\mu,\infty))}}$ is invertible  so we can cancel it and deduce  the result.
\end{pf}

We can now prove the first part of Theorem \ref{conjdt}.
By Theorem \ref{deep}
  \[\one{\P}=\one{\semi(\infty)}=\exp(\epsilon_\infty)\]
 where  $\epsilon_\infty \in (\LL-1)^{-1}
 \RHreg(\CC)_\Phi$.
The other elements in the identity of Proposition \ref{toi} lie in
$\RHreg(\CC)_\Phi$. Thus, using  Corollary \ref{ad}, we have an
equation in $\RHsc(\CC)_\Phi$
\[\H_{\leq 1} =  \H_0 \cdot \exp(\{\eta_\infty,-\}) (\Hs_{\leq 1}).\]
Applying the integration map, the Poisson bracket vanishes, so we
obtain
\[ I_\Phi(\H_{\leq 1}) = I_\Phi(\H_0) \cdot I_\Phi(\Hs_{\leq 1}),\]
and the result  follows from  \eqref{john1}, \eqref{john2} and
\eqref{newdtwo}.


\section{Rationality}
In this section we complete the proof of Theorem \ref{conjdt} by proving the rationality statement of part (ii).

\label{hardproof}


\subsection{DT invariants counting semistables}

Recall the elements $\eta_\mu$ of Theorem \ref{deep} and decompose them as
\[\eta_\mu =\sum_{\mu(\gamma)=\mu}
\eta_\gamma.\]
Following \cite[Section 6.4]{JS}, we can define invariants $N_{(\beta,n)}\in\mathbb{Q}$ counting
semistable sheaves of Chern character $(\beta,n)\in \Delta$ by
setting
\[I(\eta_{(\beta,n)})=-N_{(\beta,n)}x^{(\beta,n)}.\]
 Most of the following result was obtained previously by Joyce
and Song \cite[Theorem 6.16]{JS}; since the proof is easy we include it here.

\begin{lemma}
\label{satisfy}
Assume $\beta\in N_1(M)$ is effective and nonzero. Then the invariants $N_{(\beta,n)}$ satisfy
\begin{itemize}
\item[(a)] $N_{(\beta,n)}=N_{(\beta,n+\beta\cdot H)},$
 \item[(b)] $N_{(\beta,n)}=N_{(\beta,-n)},$
\end{itemize}
and are independent of the choice of ample divisor $H$.
\end{lemma}

\begin{pf}
Part (a) holds because tensoring with $\OO_M(H)$ defines an
automorphism of the stack $\CC$ taking semistable sheaves of class
$(\beta,n)$ to semistable sheaves of class $(\beta,n+\beta\cdot
H)$. Part (b) follows in the same way using the dualizing functor
$\DD$ and Lemma \ref{new}. Part (c) follows from Proposition
\ref{grinah} below, since the subcategory
 $\semi(I)$ corresponding to the interval
 $I=\R_{>0}\subset \R$ consists of those objects of $\Qlast$ whose quotients all have positive
Euler characteristic, and is therefore independent of the choice of the ample divisor $H$.
\end{pf}

Given a $\Delta$-graded $\kvar{\C}[\LL^{-1}]$-algebra \[A=\bigoplus_{\gamma\in \Delta} A_\gamma\]
we define an automorphism $\chi$ of $A$,
whose action we write as $a\mapsto a^\chi$, by the rule
\[a^\chi=\LL^{\chi(\gamma)}\cdot a \text{ for }a\in A_\gamma.\]
Here $\chi(\gamma)$ is the Euler characteristic of a sheaf of
class $\gamma$.

 We can consider the automorphism $\chi$ as
defining an action of $\Z$ on $A$ and form the skew-group algebra
$\Hat{A}=A*\Z$.
 In terms of generators and relations
 \[\Hat{A}=A\langle y,y^{-1} \rangle/(y * a -
 a^\chi  * y).\] Suppose
the quotient algebra
\[A_{\sc}=A/(\LL-1)A\]
is commutative, and let $\{-,-\}$ be the induced Poisson bracket.
Then the quotient
\[\Hat{A}_{\sc}=\Hat{A}/(\LL-1)\Hat{A} \] is a commutative algebra, isomorphic to $A_{\sc}[y,y^{-1}]$,
 and the induced Poisson bracket satisfies
\begin{equation}
\label{newone} \{y,a\}=\chi(\gamma)\cdot  ya\text{ for }a\in
A_\gamma.\end{equation} All this extends in the obvious way to the
topologically graded algebras $A_\Phi$.

\begin{prop}
\label{grinah} Given an interval $I\subset (-\infty,+\infty]$
bounded below, the element \[\one{\semi(I)}^{-1}*
\one{\semi(I)}^\chi\in \RH(\CC)_\Phi\] is regular and
\[I_\Phi\big(\one{\semi(I)}^{-1}* \one{\semi(I)}^\chi\big)
= \exp\bigg(\sum_{\stackrel{(\beta,n)\in\Delta}{\mu(\beta,n)\in I}}
 n N_{(\beta,n)}x^{(\beta,n)}\bigg)\in\C[\Delta]_\Phi.\]
\end{prop}

\begin{pf}
By definition of the product in the algebra
$\Hat{\RH}(\CC)_{\Phi}$
 one has
\[\add^{-1}_{\one{\semi(I)}} (y)=\one{\semi(I)}^{-1} * y * \one{\semi(I)}=\one{\semi(I)}^{-1}* \one{\semi(I)}^\chi* y. \]
The first claim then follows by the argument of Corollary
\ref{ad}.  Theorem \ref{HN} and Corollary \ref{ad} imply that
\[\add^{-1}_{\one{\semi(I)}} = \prod_{\mu\in I} \add^{-1}_{\one{\semi(\mu)}}=\prod_{\mu\in I} \exp \{-\eta_{\mu},-\}
\colon \Hat{\RH}_{\sc}(\CC)_\Phi \lra \Hat{\RH}_{\sc}(\CC)_\Phi,\]
where the product of endomorphisms are taken in order of ascending
slope, and the infinite products are interpreted as in Lemma
\ref{HN}. The homomorphism of topologically $\Delta$-graded Poisson algebras
$I_\Phi$ induces a Poisson algebra map
\[I_\Phi \colon \Hat{\RH}_{\sc}(\CC)_{\Phi} \lra \Hat{\C}[\Delta]_\Phi.\]
Applying this gives
\begin{equation}
\label{hamp} I_\Phi\big(\one{\semi(I)}^{-1}*
\one{\semi(I)}^\chi\big)\cdot y=\prod_{\mu\in I} \exp
\{-I_\Phi(\eta_{\mu}),-\} (y).\end{equation} By Lemma \ref{ady}
and \eqref{newone} this gives
\[I_\Phi\big(\one{\semi(I)}^{-1}* \one{\semi(I)}^\chi\big)\cdot y=
  \prod_{\stackrel{(\beta,n)\in\Delta}{\mu(\beta,n)\in I}} \exp
\big(n N_{(\beta,n)}x^{(\beta,n)}\big) \cdot  y\] and the result
follows.
\end{pf}

We used the following easy result.

\begin{lemma}
\label{ady}
Suppose $A$ is a topological Poisson algebra and $a,b \in A$ satisfy
\[\{a,b\}=nab.\]
 Then, providing the relevant sums converge,
\[\exp \{a,-\} (b) = \exp(na) \cdot b.\]
Furthermore, if $n=0$  then the
automorphisms $\exp \{a,-\}$ and $\exp \{b,-\}$ commute, and
\[\exp\{a+b,-\} = \exp \{a,-\} \circ \exp \{b,-\}.\]
\end{lemma}

\begin{pf}
This is left to the reader.
\end{pf}


\subsection{Rationality}

Recall the definition of the Laurent series
\[\PT_\beta(q)=\sum_{n\in\Z} \PT(\beta,n) q^n\]
 from the introduction.
To prove  Theorem \ref{conjdt}(b) we follow a strategy of Toda and prove
more, namely

\begin{thm}
\label{toda} There is an identity in $\C[\Delta]_\Phi$
\[\sum_{\beta\geq 0}\PT_\beta(-q) x^\beta=
 \exp\bigg(\sum_{\stackrel{\beta\geq 0}{n\geq 0}} n N_{(\beta,n)}x^\beta q^n \bigg)\cdot \sum_{\beta\geq 0} L_\beta(q) x^\beta,\]
where for each effective class $\beta\in N_1(M)$, the expression $L_\beta(q)$ is a Laurent polynomial in $q$ invariant under
$q\leftrightarrow q^{-1}$.
\end{thm}

This is enough by Lemma \ref{satisfy} and the following simple result.

\begin{lemma}
Fix a positive integer $d\geq 1$. Suppose given rational numbers $a_n\in\mathbb{Q}$ such that for all $n\in\Z$ one has
\[a_{n}=a_{-n}, \quad a_{n+d}=a_n.\]
Then the power series
\[\sum_{n\geq 0} n a_n q^n\]
is the Taylor expansion of a rational function in $q$ invariant
under $q\leftrightarrow q^{-1}$.
\end{lemma}

\begin{pf} This is a simple calculation; see \cite[Lemma 4.6]{To1}.\end{pf}

 Given $\mu>0$ let us
consider the following tautology in ${\RH(\CC)_\Phi}$
\begin{equation}
\label{blue}
\one{\semi([0,\mu])}^{-1} * \Hs_{\leq 1} * \one{\semi([0,\mu])}  =
\bigg(\one{\semi([0,\mu])}^{-1}* \one{\semi([0,\mu])}^\chi \bigg) * \G_\mu\end{equation}
where \[\G_\mu=
(\one{\semi([0,\mu])}^\chi)^{-1} * \Hs_{\leq 1}
 * \one{\semi([0,\mu])} .\]
Note that by Lemma  \ref{HN} and Corollary \ref{ad} the
left-hand side of \eqref{blue} is regular. Since the Poisson
bracket on $\C[\Delta]_\Phi$ is zero, it follows from
\eqref{john2} that
\begin{equation}
\label{burbage}
I_\Phi\big(\one{\semi([0,\mu])}^{-1} * \Hs_{\leq 1} *
\one{\semi([0,\mu])}\big) = I_\Phi(\Hs_{\leq 1})=\sum_{\beta\geq 0}\PT_\beta(-q) x^\beta.\end{equation}
By Proposition \ref{grinah} the first term on the right hand side is also regular, and
\begin{equation}
\label{gardoms}I_\Phi\bigg(\one{\semi([0,\mu])}^{-1}* \one{\semi([0,\mu])}^\chi \bigg)
=\exp\bigg(\sum_{\stackrel{(\beta,n)\in\Delta}{\mu(\beta,n)\in [0,\mu]}} n N_{(\beta,n)}x^{(\beta,n)}\bigg).\end{equation}
As $\mu\to \infty$ this tends to the first term in Toda's identity.
 The following result therefore completes the proof of Theorem \ref{conjdt}.

\begin{prop}
\label{barrow} The element $\G_\mu$ is regular and
\[\lim_{\mu\to \infty} I_\Phi(\G_\mu)=\sum_{\beta\geq 0} L_\beta(q) x^\beta,\]
where each $L_\beta(q)$ is a Laurent polynomial, invariant under $q\leftrightarrow q^{-1}$.
\end{prop}

\begin{pf}
Recall the subalgebra $\RH(\Qlast)\subset \RH(\CC)$, and its anti-involution $\DD$ defined in Section \ref{dualt}. Consider  the composition
\[R=\chi\circ\DD.\]
Since $\chi(\DD(E))=-\chi(E)$ for any sheaf $E\in\Qlast$, it
follows that $R$ is also an anti-involution. There is an obvious
involution of $\C[\Delta]$, which we also denote by $R$, defined
by \[R(x^{(\beta,n)})=x^{(\beta,-n)}.\] Note that by \eqref{burb},
for any $a\in \RH(\Qlast)$ one has
\begin{equation}
\label{l}I(R(a))=R(I(a)).\end{equation}

 We claim that the involution $R$ fixes the element
\[\O{\semi([-\mu,\mu])}\in \RH(\Qlast)_\Phi.\] In fact,
 for all $\gamma\in \Delta$
\begin{equation}
\label{oxf}
R(\pi_\gamma(\O{\semi([-\mu,\mu])}))=\pi_{R(\gamma)}(\O{\semi([-\mu,\mu])})\in \RH(\Qlast).\end{equation}
To see this note first
 that if $E\in \Qlast$ has class
$\gamma=(\beta,n)$ then by  Serre duality \[H^0(M,\DD(E))\isom
H^1(M,E)^*, \] and so \begin{equation} \label{verylast}
n=\chi(E)=\dim H^0(M,E)-\dim H^0(M,\DD(E)).\end{equation} Take a
finite stratification
\[\semi([-\mu,\mu])_\gamma=\big.\coprod_{r\geq 0}\semi([-\mu,\mu])_\gamma
\cap\M_r\] as in Lemma \ref{wet}. 
It follows from Lemma \ref{new}
and \eqref{verylast} that
\[\DD( \semi([-\mu,\mu])_{(\beta,n)} \cap\M_r)=\semi([-\mu,\mu])_{(\beta,-n)}\cap\M_{r-n}.\]
Applying  Lemma \ref{wet},
\[\DD( \O{\semi([-\mu,\mu])_{(\beta,n)} \cap\M_r})=\LL^n\cdot \O{\semi([-\mu,\mu])_{(\beta,-n)}\cap\M_{r-n}},\]
and the claim follows.

Next define an element
\[\J_\mu  = \Hs_{\leq 1} * \one{\semi([-\mu,\mu])}\in \RH(\CC)_\Phi.\] 
By the argument we used to prove \eqref{right}\footnote{This statement is false. See Section \ref{for2}.}
\begin{equation}
\label{ord}
\J_\mu-\O{\semi([-\mu,\mu])} \to 0 \text{ as }\mu\to \infty.\end{equation}
 Lemma \ref{HN} implies that there is an identity
\[\one{\semi([-\mu,\mu])}= \one{\semi([0,\mu])} * \one{\semi([-\mu,0))}\]
in $\RH(\CC)_\Phi$. Thus
\[\G_\mu=(\one{\semi([0,\mu])}^\chi)^{-1} *
\J_\mu
* \one{\semi([-\mu,0))}^{-1}.\]
Let us also consider the element
\[\G'_\mu=(\one{\semi((0,\mu])}^\chi)^{-1} *
\J_\mu
* \one{\semi([-\mu,0])}^{-1}.\]
  We can modify equation \eqref{blue} by
replacing the closed intervals by intervals open at 0. This has
the effect of replacing the element $\G_\mu$ by $\G'_\mu$.
Equations \eqref{burbage} and \eqref{gardoms} are unchanged by
this modification, which shows that
\[I_\Phi(\G_\mu)=I_\Phi(\G'_\mu).\]

Using Lemma \ref{new} and the fact that $R$ is an
anti-automorphism we find
\[R(\G'_\mu)=(\one{\semi([0,\mu])}^\chi)^{-1} *
R(\J_\mu)
* \one{\semi([-\mu,0))}^{-1}.\]
It follows from \eqref{oxf} and \eqref{ord} that
\[\G_\mu-R(\G'_\mu)=  (\one{\semi([0,\mu])}^\chi)^{-1} *
(\J_\mu-R(\J_\mu))
* \one{\semi([-\mu,0))}^{-1}\to 0\]
as $\mu\to\infty$. Thus by \eqref{l}
\[I_\Phi(\G_\mu)-R(I_\Phi(\G_\mu))\to 0\]
On the other hand, the limit of $I_\Phi(\G_\mu)$ as $\mu\to \infty$ certainly exists because the other
factors in \eqref{blue} have convergent integrals. Thus we conclude that this limit is invariant under $R$ and the result follows.
\end{pf}

%
%
%

\subsection{Correction added after publication}

\label{for}
Yukinobu Toda kindly pointed out a mistake in the published version of this paper. Fortunately it is easily fixed. The mistake lies in the statement of Theorem 6.3. This is only true as stated if we change the definition of the regular Hall algebra in Section 5.1. Rather than defining this as the span of the given symbols over $\kvar{\C}[\LL^{-1}]$, we should take the span over the larger subring\begin{equation}
\label{whoops}\kvar{\C}[\LL^{-1}][ (1+\LL+\cdots + \LL^k)^{-1} : k\geq 1]\subset \kst{\C}.\end{equation} That is, we should invert the classes of all projective spaces. The definition of the semi-classical Hall algebra $\RHsc(\CC)$ is unaffected by this change, and everything else goes through as before.

To illustrate the problem, consider the following example which was communicated to the author by Toda.  Suppose there is a single stable object $E$ of slope $\mu$ which is moreover spherical, so that $\Ext^1_X(E,E)=0$. For example, one could take $E=\OO_C$, where $C\subset X$ a rational curve with normal bundle $\OO_C(-1)^{\oplus 2}$. Every semistable object of slope $\mu$ is  then isomorphic to $E^{\oplus n}$ for some $n\geq 0$ so
\[1_{\semi}(\mu) = \bigsqcup_{n\geq 0} [X_n\subset \M],\]
where $X_n\subset \M$ is the open  substack of objects isomorphic to $E^{\oplus n}$. Clearly $X_n$ is the quotient stack $\Spec(\C)/\GL(n,\C)$. Since $\Hom_X(E,E)=\C$, the fibre of the map $(a_1,a_2)$ in the diagram \eqref{harrogate} over the point $(E,E)$ is 
the quotient stack $\Spec(\C)/\C$. Thus if we pull back the Hall algebra element $\epsilon_\mu$ via the map $\Spec(\C)\to \M$ corresponding to the object $E^{\oplus 2}$, we get a contribution
\[\frac{1}{[\GL(2,\C)]}-\frac{1}{\LL\cdot [\GL(1,\C)]^2}= \frac{1}{\LL-1}\cdot \frac{-1}{2\,\LL\,(\LL+1)}.\]
Multiplying by $\LL-1$ we get an element $\eta_\mu$ which is not regular in the sense used in the published version, but which is regular in our new sense, since it lies in the ring \eqref{whoops}.

The argument we gave in the proof of Theorem 6.3  allows one to conclude that the  element $\eta_\mu\in \H(\CC)$ is a $\Lambda$-linear combination of maps from schemes, where $\Lambda=\kst{\C}$. But this does not imply Theorem 6.3 in either the  published or modified versions. However applying the framework of \cite[Section 6]{Jo0} with the subring $\Lambda^\circ\subset \Lambda=\kst{\C}$ given by \eqref{whoops} allows us to conclude that $\eta_\mu$  is a $\Lambda^\circ$-linear combination of maps from schemes, which is precisely what our modified version of Theorem 6.3 claims.

\subsection{A second correction}
\label{for2}
Yukinobu Toda communicated a second error in the published version of this paper, which was found by an undergraduate(!) student, Tasuki Kinjo. The problem is that the argument used to prove \eqref{right} does not extend to give \eqref{ord} as claimed, since the Euler chararacteristic $\chi(A)$ of a stable pair $\OO_X\to A$ with fixed curve class $\ch_2(A)=\beta$ is not bounded above. In fact, a simple example shows that  equation \eqref{ord} is false in general. 

To fix the proof of Proposition \ref{barrow} we require some additional arguments which we now explain.
We first follow the argument as given above, ignoring the false claim \eqref{ord}, until we get to the line
\begin{equation}
\label{line}\G_\mu-R(\G'_\mu)=  (\one{\semi([0,\mu])}^\chi)^{-1} *
(\J_\mu-R(\J_\mu))
* \one{\semi([-\mu,0))}^{-1}.\end{equation}
To conclude we must show that this expression tends to zero as $\mu\to \infty$. Let us fix a pair $(\beta,n)$. We must show that  the projection $\pi_{(\beta,n)}$ applied to \eqref{line} is zero for sufficiently large $\mu>0$.

The first observation is that  the projections $\pi_{(\beta,n)}$ applied to the left-hand side of \eqref{line} is independent of $\mu$ for large enough $\mu>0$. This follows immediately from the definitions of the elements $\G_\mu$  and $\G'_\mu$.  Let us then fix some such sufficiently large $\mu>0$. Since $\pi_{(\beta,n)} (\G_\mu)$ is of finite type we can then take $s>0$ large enough so that the projection of the left-hand side of \eqref{line}  is represented by a map into the open substack $\semi([-s,s])\cap \M_{(\beta,n)}\subset \M_{(\beta,n)}$.

We now  apply  the the argument used to prove \eqref{right}. Increasing $\mu$ if necessary, it shows that for  any $0< \beta'<\beta$, and any $n'\in \Z$, the projection $\pi_{(\beta',n')} (\J_\mu - R(\J_\mu))$ becomes zero when restricted to the open substack $\semi([-\mu,\mu])$. Indeed, over this open substack  one has $\pi_{(\beta',n')}(\J_\mu)= \pi_{(\beta',n')}(1^\OO)$.

To calculate the product \eqref{line} we must consider a filtration of the form \[0=E_0\subset E_1\subset E_2 \subset E_3=E,\]\[ E_1/E_0 = A \in \semi([0,\mu]), \qquad E_2/E_1 = B, \qquad E_3/E_2 = C \in \semi([-\mu,0)).\]
The claim is that for $\mu\gg s$ we have $E \in \semi([-s,s])$ implies $B \in \semi([-\mu,\mu])$.  It follows immediately from this claim and the second observation above that the product \eqref{line} is zero when restricted to the substack $\semi([-s,s])$. But by the first observation above the product is then necessarily zero.

To prove the claim,  we can assume that  $\mu>rs$ where $r = \ch_2(\beta)\cdot H$ (recall the definition of $\mu$ from Section \ref{ex}). Suppose for a contradiction that there is a subobject $K\subset B$ of slope $>\mu$. Its inverse image in $E_2$ is a subobject $M$ which satisfies $\chi(M)\geq\chi(K)$, and hence  $\mu(M)>\mu(K)/r$.  But now $M\subset E$, so $\mu(M) >\mu(K)/r > \mu/r$ which contradicts $E \in \semi([-s,s])$ if $\mu > rs$. A similar argument applies when $B$ has a quotient of slope $<-\mu$.

\end{document}